\documentclass[12pt]{amsart}
\usepackage{mathrsfs}
\usepackage{amsmath,amssymb,amsfonts,latexsym,txfonts}
\usepackage{eucal,pifont}
\usepackage{fancyhdr,tikz}

\newlength{\originalbase}
\input amssym.def
 
\input amssym.tex

\renewcommand{\SS}{\mathbb{S}}
\newcommand{\RR}{\mathbb{R}}

\newcommand{\finial}{$$\mbox{\ding{167}}$$}

\begin{document}

\begin{center}
{\Large\bf A probabilistic proof of the spherical excess formula}\\[6mm]
{\bf\em 
Daniel A.~Klain}
\end{center}

\vspace{4mm}

A triangle $T$ in the unit sphere with 
inner angles $\theta_1$, $\theta_2$, and $\theta_3$ 
has area given by the 
{\em spherical excess formula:}\footnote{This formula was discovered in 1603 by Thomas Harriot \cite[p. 65]{Still} and is also known as Girard's formula \cite[p. 95]{Cox-geom}.}  
%A proof using the inclusion-exclusion principle is given in \cite[p. 158]{Lincee}.} 
\begin{align}
\mathtt{Area}(T) \;=\; \theta_1+\theta_2+\theta_3-\pi.
\label{exform}
\end{align}
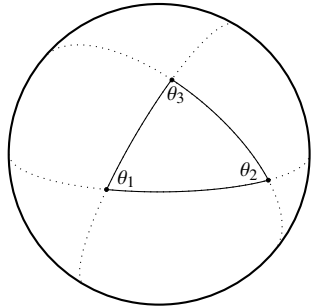
\begin{figure}[ht]
\centering
\begin{tikzpicture}[xscale=1,yscale=1]
\draw[thick] (0,0) circle [radius=2];
\draw[dotted,domain=180:360] plot ({2*cos(\x)}, {0.5*sin(\x)});
\draw[dotted,domain=0:180] plot ({cos(60)*2*cos(\x) - sin(60)*0.4*sin(\x)}, {sin(60)*2*cos(\x) + cos(60)*0.4*sin(\x)});
\draw[dotted,domain=0:180] plot ({cos(-40)*2*cos(\x) - sin(-40)*0.9*sin(\x)}, {sin(-40)*2*cos(\x) + cos(-40)*0.9*sin(\x)});
\draw[smooth,domain=249.5:316.5] plot ({2*cos(\x)}, {0.5*sin(\x)});
\draw[smooth,domain=62:112] plot ({cos(60)*2*cos(\x) - sin(60)*0.4*sin(\x)}, {sin(60)*2*cos(\x) + cos(60)*0.4*sin(\x)});
\draw[smooth,domain=48:105] plot ({cos(-40)*2*cos(\x) - sin(-40)*0.9*sin(\x)}, {sin(-40)*2*cos(\x) + cos(-40)*0.9*sin(\x)});
\draw[fill] (-0.7,-0.4683) circle [radius=0.025];
\draw[fill] (1.4507,-0.3442) circle [radius=0.025];
\draw[fill] (0.17,0.9897) circle [radius=0.025];
\node at (-0.7+.27,-0.4683+.15) {\tiny $\theta_1$};
\node at (1.4507-.25,-0.3442+.12) {\tiny $\theta_2$};
\node at (0.17+.06,0.9897-.2) {\tiny $\theta_3$};
\end{tikzpicture} 
\caption{A spherical triangle.} 
\label{s-tri}
\end{figure}

\begin{figure}[ht]
\centering
\begin{tikzpicture}[xscale=1,yscale=1]
\shade[top color=gray, bottom color=white] (0-1.33,-2)--(0,0)--(0+0.6,-2)--(0-1.33,-2);
\draw[fill] (0,-0.01) circle [radius=0.03];
\draw[thick,->] (0,0)--(-1.33,-2);
\draw[thick,->] (0,0)--(-0.4,-2.5);
\draw[thick,->] (0,0)--(0.6,-2);
\draw[dashed] (-0.99,-1.5)--(-0.3,-1.875)--(0.45,-1.5);
\draw[dotted] (0.45,-1.5)--(-0.99,-1.5);
\draw[ultra thick,dashed,->] (1.6,-1)--(3,-1);
\node [above] at (2.25,-0.8) {\bf ?};
\shade[top color=gray, bottom color=white] (5-1.33,-2)--(5,0)--(5+0.6,-2)--(5-1.33,-2);
\draw[fill] (5,-0.01) circle [radius=0.03];
\draw[thick,->] (5,0)--(5-1.33,-2);
\draw[thick,dotted,->] (5,0)--(5-0.4,-2.5);
\draw[thick,->] (5,0)--(5+0.6,-2);
\node at (4.7,-2.9) {\scriptsize 1 vertex};
\node at (4.7,-3.23) {\scriptsize 2 edges};
\shade[outer color=white, inner color=gray] (7.5,-1) circle [radius=1]; 
\draw[fill, color=white] (6.5,-1)--(8.5,-1)--(8.5,1)--(6.5,1)--(6.5,-1);
\draw[fill] (7.5,-1) circle [radius=0.03];
\draw[thick,->] (7.5,-1)--(8.5,-1);
\draw[thick,->] (7.5,-1)--(6.5,-1);
\draw[thick,dotted,->] (7.5,-1)--(7.2,-2);
\node at (7.5,-2.9) {\scriptsize (measure zero};
\node at (7.5,-3.23) {\scriptsize outcome)};
\shade[outer color=white, inner color=gray] (10,-1) circle [radius=0.7];
\draw[fill] (10,-1) circle [radius=0.03];
\draw[thick,dotted,->] (10,-1)--(9.4,-0.5);
\draw[thick,dotted,->] (10,-1)--(10.8,-0.7);
\draw[thick,dotted,->] (10,-1)--(10,-1.7);
\node at (10,-2.9) {\scriptsize 0 vertices};
\node at (10,-3.23) {\scriptsize 0 edges};
\end{tikzpicture} 
\caption{Projections of a 3-dimensional cone.} 
\label{cone-proj}
\end{figure}
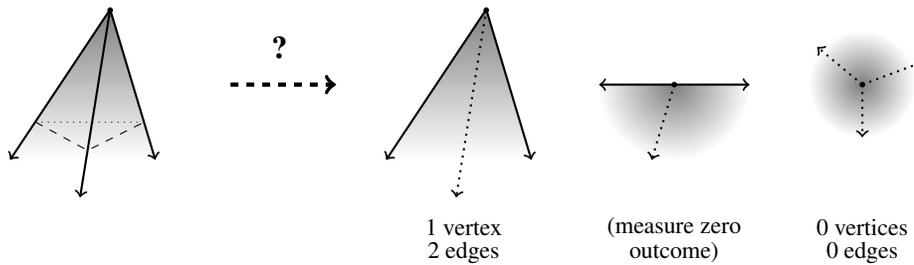

This note offers a probabilistic proof of the angle excess formula~(\ref{exform}), 
based on the observation that an unbounded %triangular 
cone at the origin in $\RR^3$ 
has only three kinds of $2$-dimensional orthogonal projections: 
a cone in $\RR^2$, a half-plane in $\RR^2$ (this is an outcome of probability zero), and all of $\RR^2$. See Figure~\ref{cone-proj}.

Observe that, if we omit the middle outcome of measure zero, the number of edges on each projected figure is {\em twice} the number of vertices.
%\finial

Some notation will help to interpret angles as probabilities.  
Let $\SS$ denote the unit sphere in $\RR^3$ centered at the origin, having surface area $4\pi$. 

Suppose that $P$ is a convex polytope in $\RR^3$, and let $x$ be any point of $P$.
The {\em solid inner angle} $a_P(x)$ of $P$ at $x$ is given by
$$a_P(x) = \{ u \in \SS \; | \; \ x+\epsilon u \in P \;\; \hbox{ for some } \epsilon > 0 \}.
$$
Let $\alpha_P(x)$ denote the measure of the solid angle $a_P(x) \subseteq \SS$, 
given by the usual surface area measure on subsets of the sphere.

If $F$ is a proper face of a convex polytope $P$,
then the solid inner angle measure $\alpha_P(x)$
is the same at every point $x$ in the relative interior of $F$. 
This value will be denoted by
$\alpha_P(F)$.

Consider 
the case of an unbounded cone $C$ with single vertex at the origin $o$, as in Figure~\ref{cone-proj}.
Specifically, let $v_1, v_2, v_3$ be three linearly independent unit vectors in $\RR^3$, 
and let $C$ denote all non-negative linear combinations:
$$C = \{t_1 v_1 + t_2 v_2 + t_3 v_3 \;|\; t_i \geq 0\}.$$
The polyhedral cone $C$ has exactly one vertex at $o$ 
and three (unbounded) edges $e_i$ along the directions of the vectors $v_i$.
Note that $\alpha_C(o)$ is the area of the spherical triangle with vertices at $v_i$.  Denote the inner angles 
of this triangle by $\theta_i$, 
as in Figure~\ref{s-tri} (where $o$ lies at the center of the sphere in Figure~\ref{s-tri}).

Given a uniformly random unit vector $u$, let $C_u$ denote the orthogonal projection of $C$ onto the plane $u^\perp$.  Evidently $C_u$
will resemble one of the outcomes in Figure~\ref{cone-proj}.  Specifically, $C_u$ will cover the entire plane $u^\perp$ iff $u$ lies in the
interior of $\pm a_C(o)$.  It follows that $C_u = u^\perp$ with probability 
$$\frac{\mathtt{Area}(a_C(o))+\mathtt{Area}(-a_C(o))}{4\pi} \;=\; 
\frac{2\alpha_C(o)}{4\pi} \;=\; \frac{\alpha_C(o)}{2\pi}.$$
Since the number of vertices of $C_u$ is either $0$ or $1$, the expected number of vertices of $C_u$ is given by
the complementary probability
\begin{align}
E(\hbox{\# of vertices}) \;=\; 1- \frac{\alpha_C(o)}{2\pi}.
\label{expv}
\end{align}

%\finial

Meanwhile, an edge $e$ projects to the interior of $C_u$ iff $u$ lies in 
the interior of $\pm a_C(e)$.  Taking the complement as before, $e$ projects to a boundary edge of $C_u$
with probability $1- \frac{\alpha_C(e)}{2\pi}$.
Observe that each solid angle measure $\alpha_C(e_i)$ is given by $2\theta_i$ (see Figure~\ref{lune}),
so that the expected number of edges of $C_u$ is 
\begin{align}
E(\hbox{\# of edges}) \;=\;
\sum_{i} \left( 1 - \frac{\alpha_C(e_i)}{2\pi} \right)
\;=\;
\sum_{i} \left( 1 - \frac{\theta_i}{\pi} \right)
\;=\;
3 - \frac{\theta_1+\theta_2+\theta_3}{\pi}.
\label{expe}
\end{align}
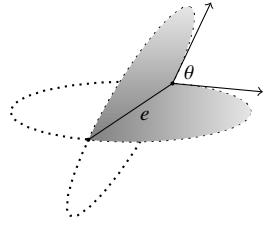
\begin{figure}[ht]
\centering
\begin{tikzpicture}[xscale=0.8,yscale=0.8]
\draw[thick,dotted,domain=0:360] plot ({2*cos(\x)}, {0.5*sin(\x)});
\draw[thick,dotted,domain=0:360] plot ({cos(60)*2*cos(\x) - sin(60)*0.4*sin(\x)}, {sin(60)*2*cos(\x) + cos(60)*0.4*sin(\x)});
\shade[left color={rgb:black,5;white,5}, right color={rgb:black,1;white,9},domain=-110.5:69.5,variable=\x] 
(-0.7,-0.4683) -- plot ({2*cos(\x)}, {0.5*sin(\x)})
-- (0.7,0.4683) -- cycle;
\shade[top color={rgb:black,1;white,9}, bottom color={rgb:black,5;white,5},domain=-68:112,variable=\x] 
(-0.7,-0.4683) -- plot ({cos(60)*2*cos(\x) - sin(60)*0.4*sin(\x)}, {sin(60)*2*cos(\x) + cos(60)*0.4*sin(\x)})
-- (0.7,0.4683) -- cycle;
\draw[thin] (-0.7,-0.4683)--(0.7,0.4683);
\draw[fill] (-0.7,-0.4683) circle [radius=0.025];
\draw[fill] (0.7,0.4683) circle [radius=0.025];
\node at (0.7+.27,0.4683+.19) {\tiny $\theta$};
\node at (-0.06+.3,-0.06+0) {\tiny $e$};
\draw[->] (0.7,0.4683)--(0.7+0.7974*.8,0.4683+1.681*.8);
\draw[->] (0.7,0.4683)--(0.7+1.873*.8,0.4683-0.175*.8);
\end{tikzpicture} 
\caption{$\alpha_C(e)= 2\theta$} 
\label{lune}
\end{figure}

Since the number of edges in $C_u$ is almost surely {\em twice} the number of vertices
(see Figure~\ref{cone-proj}), the identities~(\ref{expv}) and~(\ref{expe}) imply that
\begin{align}
3 - \frac{\theta_1+\theta_2+\theta_3}{\pi} 
\;=\; E(\hbox{\# of edges})
\;=\; 2E(\hbox{\# of vertices})
\;=\; 2- \frac{\alpha_C(o)}{\pi}.
\label{almost}
\end{align}

It is now immediate from~(\ref{almost}) that
$$\alpha_C(o) \;=\; \theta_1 + \theta_2 + \theta_3 - \pi,$$
as asserted in~(\ref{exform}).
\finial

While the argument above generalizes easily to the case of a spherical polygon with more than 3 sides,
in higher dimensions there is a proliferation of cases that makes this approach much more complicated.
However, a variation of this approach
was applied in~\cite{F-K} to an $n$-simplex $\Delta$ to obtain the identity
\begin{align*}
p_\Delta \;=\; \frac{2}{n\omega_n} \sum_v \alpha_\Delta(v),
%\label{anglesum}
\end{align*}
where $p_\Delta$ is the probability that a random orthogonal projection $\Delta_u$ is an $(n-1)$-simplex,
$\omega_n$ denotes the volume of the unit ball in $\RR^n$, 
and where the sum is taken over all vertices $v$ of the $n$-simplex $\Delta$.  

Similar and more general arguments were used even earlier by Perles and Shephard \cite{Shep-Perles}
(see also \cite[p.~315a]{Grunbaum} and \cite{Welzl})
to give a simple proof of the Gram-Euler identity for convex polytopes,
\begin{align*}
\sum_{F \subseteq \partial P} (-1)^{\dim F} \alpha_P(F) \;=\; (-1)^{n-1} n \omega_n,
\end{align*}
where the sum is taken over all proper faces $F$ of an
$n$-dimensional convex polytope $P$.

\vspace{2mm}
\noindent {\em Department of Mathematical Sciences, 
University of Massachusetts Lowell, Lowell, MA 01854 USA \\
Daniel\_Klain@uml.edu}

\end{document}